\documentclass[11pt]{article}

\usepackage{amsmath, graphicx, latexsym, amssymb}
\usepackage{epsfig}
\usepackage{float}

\newtheorem{Df}{Definition}[section]

\newtheorem{Th}[Df]{Theorem}
\newtheorem{Co}[Df]{Corollary}

\newtheorem{Pro}[Df]{Proposition}

\newtheorem{Ex}[Df]{Example}
\newtheorem{Exs}[Df]{Examples}

\newtheorem{Rem}[Df]{Remark}

 \newsavebox{\no} \savebox{\no}[.18in]{$=\!\!\!\!\! /$}
 \newsavebox{\lno}
\savebox{\lno}[.13in]{ $\! \mbox{\scriptsize =}\!\!\!
\mbox{\scriptsize  /}~$} 
\savebox{\nexists}[.18in]{$\exists\!\!\!\! /$}

\newcommand{\proof}{\noindent{\bf Proof:~}}

\newcommand{\skp}{\vspace{\baselineskip}}
\newcommand{\noi}{\noindent}

\usepackage{amssymb}
\usepackage{amsmath}
\usepackage{graphicx}
\usepackage{wrapfig}
\usepackage{bbm}
\usepackage{float}
\usepackage{tikz}
\usetikzlibrary{arrows,calc,shapes,decorations.pathreplacing}
\usepackage{wrapfig}
\usepackage{hyperref}
\hypersetup{pdftex,colorlinks=true,allcolors=blue}
\usepackage{hypcap}
\usepackage{caption}
\usepackage{pgfplots}
\usepackage{subcaption}

\begin{document}

\baselineskip 20pt
\title{On the price of anarchy in a single server queue with heterogenous service valuations induced by travel costs}
\author{
Refael Hassin\\ Department of Statistics and Operations Research, \\Tel Aviv
University, \\ Tel Aviv 69987, Israel.
\and Irit Nowik \ \ \ \ Yair Y. Shaki \footnote{hassin@post.tau.ac.il, nowik@jct.ac.il y\_shaki@hotmail.com }\\Department of Industrial Engineering and Management, \\The Jerusalem College of Technology, \\Jerusalem, Israel.
 }


%


\maketitle

\noi{\bf Abstract}

This work presents a variation of Naor's strategic observable model (1969), by adding a component of customer heterogeneity induced by the location of customers in relation to the server.
Accordingly, customers incur a ``travel cost'' which depends linearly on the distance of the customer from the server.
The arrival of customers with distances less than $x$ is assumed to be a Poisson process with rate $\lambda(x)=\int_0^x h(y) dy<\infty,$ where
$h(y)$ is a nonnegative  ``intensity'' function of the distance $y.$
In a loss system M/G/1/1 we define the threshold Nash equilibrium strategy $x_e,$ and the optimal social threshold strategy $x^*.$
We consider the price of anarchy (PoA) and prove that it converges to 1 when $x_e\to 0.$
The behaviour of PoA when $x_e\to\infty$ is more complex and interesting. We show that if the rate of arriving customers is bounded then PoA converges to 1, i.e., in the limit there is no difference between the social and equilibrium optimal benefits (even though the corresponding optimal strategies $x^*$ and $x_e$ do not coincide). The rest of the paper is dedicated for the case in which the rate of arriving customers is unbounded. We develop an explicit  formula to calculate $\lim\limits_{x_e\to \infty}PoA(h, x_e)$ when it exists. We present sufficient conditions for the limit to exist and for the existence of a simple formula for calculating it. We prove that if the relation between two intensity functions converges to a positive constant, then the corresponding limits of PoA coincide. If, on the other hand, one intensity function is larger than the other from some point on, then under certain conditions the limit of PoA (if exists) will be larger for the larger intensity function.
For all intensity functions $h,$
we prove that if $h$ converges to a constant then PoA converges to $2,$ and that
 if from some point on $h$ decreases (increases) monotonically then the limit of PoA, if exists, is smaller (larger) than $2.$
In a system with a queue we prove that the price of anarchy may be unbounded already in the simple case of uniform arrival, namely $h\equiv c,$ where $c>0.$

\skp

\noi{\bf Keywords:} profit maximization; price of anarchy, travel costs, observable queue.

\section{Introduction}

The performance and optimization of service systems has attracted much attention in recent years (see \cite{HH}, \cite{H1}).
Naor \cite{N} was the first to introduce a queueing model that describes customer rational decisions. The model considers an FCFS M/M/$1$ system and defines social and individual welfare.
The assumptions in Naor's model are:
\begin{enumerate}
\item A stationary Poisson stream of customers with parameter $\lambda$.
\item Service times are independent and exponentially distributed with parameter $\mu$.
\item There is a cost of $C$ per unit time spent while waiting or in service.
\item The benefit from a completed service is $R$.
\end{enumerate}

The model's parameters can be normalized so that there are only two relevant parameters:

\begin{itemize}

\item The utilization factor $\rho=\frac{\lambda}{\mu}$.

\item The reward in terms of the expected waiting cost during a service $\frac{R\mu}{C}$.

\end{itemize}

The Nash equilibrium solution in this model is very simple since there exists a dominant pure threshold strategy $n_e,$ such that an arriving customer joins the queue if and only if the observed queue upon arrival is shorter than $n_e.$
This strategy maximizes the individual's expected welfare regardless of the strategies adopted by the others. Similarly, the optimal behavior, namely the strategy that maximizes social welfare (i.e., the sum of gains) is characterized by a pure threshold strategy $n^*$.

Naor observed that the  socially optimal threshold is  bounded  by the Nash equilibrium strategy. Namely,
$$n^* \leq n_e.$$

The ``Price of Anarchy'' PoA, measures the inefficiency of selfish behaviour.
It is defined as the ratio of the social welfare under optimum to the social welfare under equilibrium.

 Naor assumes that customers are homogeneous with respect to
 service valuation. Most of the more recent works on observable queues
 (i.e., assuming customers know the queue length before deciding to join
 it) follow this assumption. Some exceptions are described in section
 2.5 of \cite{HH}. For example, Larsen \cite{Larsen}
 assumes that the service value is a continuous random variable and
 proves that the profits and social welfare are unimodal functions of
 the price. For the case of a loss system (where customers join iff
 the server is idle) Larsen proves that the profit-maximizing fee
 exceeds the socially optimal fee. Miller and Buckman \cite{MB} consider an
 M/M/$s$/$s$ loss system with heterogeneous service values and characterize the socially optimal fee.

Gilboa-Freedman, Hassin and Kerner~\cite{GHK} discuss the PoA in Naor's model. They find that it has an odd behavior in two aspects: First, it increases sharply (from 1.5 to 2) as the arrival rate comes close to the service rate; Secondly, it becomes unbounded exactly when the arrival rate is greater than the service rate, which is odd since the system is always stable.

 In this paper we introduce heterogeneity in service valuation through
 a Hotelling-type model where customers reside in a ``linear city'' and
 incur ``transportation costs'' from their locations to the location of
 the server.
 Similar models have been investigated (e.g. \cite{AH, DS, EN, H2, PS, SREMJ}) but they all assume
 a constant density (possibly restricted to an interval).

Kwasnica and Stavrulaki \cite{KE} (2008) assume that the customers are homogeneous with respect to price and 
the shipping delay for a customer at distance $s$ from the server is a linear function
$g(s) = G_0 + G \cdot s$ where $G_0 \ge 0$ and $G > 0$. They obtain a solution for the optimal capacities and service semirange in the monopolistic case.

Gallay and Hongler \cite{GO}(2008) assume two observable queues
with general service rates. Servers are located at points $x_i \in
[-\Delta,+\Delta], i=1, 2$.
Customers arrive at rate $\Lambda$ from locations uniformly distributed over
$[-\Delta,+\Delta]$.
An arriving customer located at $x$ chooses which server to join by comparing the
expected utilities $U_i(x)$ from joining the ith queue, $i = 1, 2$:
$$U_i(x) = a - p_i - c_t|x - x_i| - c_w E(W_i|N_i),$$
where $a$ is service value, $p_i$ is service fee at
server $i$, $c_t$ is transportation cost per unit distance (though travel is instantaneous), $c_w$
is waiting-cost rate, $N_i$ is queue length and $E(W_i|N_i)$ is conditional expected sojourn time at $i$.
They restrict the analysis to heavy traffic.

 In contrast
 we allow more general intensity functions. Namely, the arrival of customers with distances less than $x$ is assumed to be a Poisson process with rate $\lambda(x)=\int_0^x h(y) dy<\infty,$ where
$h(y)$ is a nonnegative  ``intensity'' function of the distance $y.$
 The intensity function and
 (w.l.o.g linear) travel costs jointly generate the distribution of
 customer service valuations. A simple example is a two-dimensional city, in which the arrival of customers is uniform. In this case the intensity function can be defined as $h(x)=2\pi x,$ and so the arrival of customers with distances less than $x$ is assumed to be a Poisson process with rate $\lambda(x)=\int_0^x 2\pi y  dy = \pi x^2.$    
 
 We find that this model of heterogeneity
 is preferred to an exogenous distribution of valuations as we
 can use natural assumptions on density for generating the distribution
 rather than assuming it is exogenous.

We first consider an M/G/1/1 loss system and define $x_e$ as the threshold Nash equilibrium strategy, namely the maximal distance from which customers will join the queue under individual optimality, and $x^*$ as the threshold value that attains optimal social welfare. We show how $x_e$ is determined by the parameters $R, \mu, c_t$ and $c_w,$ of the model.

For any nonnegative intensity function $h,$ we consider the price of anarchy (PoA) and prove that $PoA\to 1$ when $x_e\to 0.$
The behaviour of PoA when $x_e\to\infty$ is more complex and interesting. We show that this limit does not always exist, and that it may be infinite. We show that  if the amount of customers arriving from far is small (i.e., $\int_0^{\infty} h(y) dy <\infty$), then in the limit there is no difference between the social and equilibrium optimal benefits, namely $\lim\limits_{x_e\to \infty}PoA(h, x_e)=1,$ (even though the corresponding optimal strategies $x^*$ and $x_e$ do not coincide). The rest of the paper is dedicated for the case in which $\int_0^{\infty} h(y) dy =\infty.$ We develop an explicit  formula to calculate $\lim\limits_{x_e\to \infty}PoA(h, x_e)$ when it exists and show that if $h,h'$ are monotonic then this limit exists and we arrive at a very simple formula to calculate it. We show that if $h$ converges to a  constant then $\lim\limits_{x_e\to \infty}PoA(h, x_e)=2.$
We prove that if $h$ decreases (increases) monotonically and $\lim\limits_{x_e\to \infty}PoA(h, x_e)$ exists, then $\lim\limits_{x_e\to \infty}PoA(h, x_e)\leq 2$ ($\geq 2$).  In addition, for any two nonnegative intensity functions $h_1, h_2$ s.t. $h_1/h_2\to c>0,$ we prove that if the corresponding  $\lim\limits_{x_e\to \infty}PoA(h_i, x_e), \ \ i=1,2,$ exist, then $\lim\limits_{x_e\to \infty}PoA(h_1, x_e)=\lim\limits_{x_e\to \infty}PoA(h_2, x_e).$ Finally, if $h_1, h_2,h'_1, h'_2,$ are all monotonic and from some point on $h_1\leq h_2,$ then the corresponding  $\lim\limits_{x_e\to \infty}PoA(h_i, x_e), \ \ i=1,2,$ exist, and: $\lim\limits_{x_e\to \infty}PoA(h_1, x_e)\leq\lim\limits_{x_e\to \infty}PoA(h_2, x_e).$

In a system with a queue we prove that the price of anarchy may be unbounded already in the simple case of uniform arrival, namely $h\equiv c,$ where $c>0.$

\section{Model description}

Consider an M/M/$1$ queue, with a server located at the origin. 
The model makes the
following assumptions:
\begin{enumerate}
\item For all $x\geq 0,$ customers with distances less than $x,$ arrive to the system according to a Poisson process with rate $\lambda (x)=\int\limits_0^x h(y)dy,$ where $h(y)$ is an  ``intensity'' function defined for all $y\geq 0.$

\item An ``intensity'' function $h$ may be any nonnegative function for which $\lambda (x)=\int\limits_0^x h(y)dy$ is finite for all $x\geq 0.$

\item Each customer knows his distance from the
server.

\item The queue length is observable.

\item Customers are risk neutral, maximizing expected net benefit.

\item The service distribution is general with average rate $\mu$.

\item The benefit from a completed service is $R$.

\item The waiting cost is $c_w$ per unit  time (while in the system).

\item The traveling cost is $c_t$ per unit  distance and traveling is instantaneous. (If $c_t=0$, we obtain the original Naor's model with rate $\lambda=\int\limits_0^\infty h(y)dy$.)
\item $\nu=\frac{R\mu}{c_w}>1$.
\item The decision of the customer is whether to join the queue or balk.

\end{enumerate}

\begin{Rem}\label{trav}
Although this model is described as a model with locations and travel costs, it may alternatively describe a model with heterogenous services values. In this interpretation,
the value of service, to a customer of ``type $x$''  is: $R-c_tx.$
\end{Rem}

\section{A loss system}\label{Loss}

First, we consider an M/M/$1/1$ loss system i.e., if the server is
busy, customers will balk. The optimal strategy of a customer
located at a distance $x$ from the origin, is to arrive if the server is idle and $R
\ge \frac{c_w}{\mu}+c_tx.$ Equivalently, $x\le \frac{R\mu -
c_w}{c_t\mu}=\frac{\frac{R\mu}{c_w} -
1}{\frac{c_t\mu}{c_w}}=\frac{\nu-1}{\kappa},$ where $\nu=\frac{R\mu}{c_w},$ and  $\kappa=\frac{c_t\mu}{c_w}$. Consequently, the threshold strategy
\begin{equation}\label{xe}
x_e=\frac{R\mu -
c_w}{c_t\mu}=\frac{\nu-1}{\kappa}
\end{equation}
is the unique Nash equilibrium strategy. Under
this individual strategy, a customer located at a distance $x,$
enters service \emph{iff} the
server is idle and $x\le x_e$. If the server is idle then the utility of the arriving customer located at point $x,$ is:
\begin{equation}\label{EI}
R- \frac{c_w}{\mu}-c_tx=c_t(x_e-x).
\end{equation}

Define $\rho(x)$ as: $$\rho(x)=\frac{\lambda(x)}{\mu}=\frac{1}{\mu}\int\limits_0^x h(y)dy.$$
 Then the probability of an idle server is:
$$\pi_{0}(x)=\frac{1}{1+\rho(x)}= \frac{1}{1+\frac{1}{\mu}\int\limits_0^x h(y)dy}.$$
Thus by~(\ref{EI}), the expected social benefit  per unit of time associated
with threshold  $x$ satisfies the following equation:
\begin{equation}\label{socialbenefit}
S(x)=c_t\int_{0}^x \left(x_e- y\right) h(y)\pi_0(x)dy=
 \frac{c_t\int_{0}^x (x_e-y)h(y)dy }{1+\frac{1}{\mu} \int_{0}^x h(y)dy}.
\end{equation}
Let $x^*$ be the threshold value that attains optimal social welfare.

Note that between $x_e$ and $x^*,$ \ $x_e$ is the more basic feature of the model in the sense that it is determined by the parameters of the model such that   \ $x_e=\frac{R\mu -
c_w}{c_t\mu}$  \ (see~(\ref{xe})). Hence we relate to $x^*$ as a function of $x_e.$ 

Denote: $$f(x)=\frac{1}{\mu}\int_{0}^x (x-y)h(y)dy +x.$$
\begin{Pro}\label{P3000}
For every $x_e\geq 0$ the optimal threshold strategy $x^*$ is unique and satisfies:
\begin{itemize}
\item $x^*=f^{-1}(x_e),$
\item $x^*\geq 0,$ 
\item $\lim\limits_{x_e\to\infty}f^{-1}(x_e)=\infty.$
\end{itemize}
\end{Pro}

\proof
In order to find $x^*$ we compute:
$$S'(x)=\frac{ c_t\Bigl( (x_e-x)h(x)[1+\frac{1}{\mu} \int_{0}^x h(y)dy] \Bigr)- c_t\Bigl(\int_{0}^x (x_e-y)h(y)dy \Bigr)[\frac{1}{\mu} h(x)]}{[1+\frac{1}{\mu} \int_{0}^x h(y)dy]^2}=0. $$
We obtain:
$$ (x_e-x)\Bigg[1+\frac{1}{\mu}\int_{0}^x h(y)dy\Bigg] =\frac{1}{\mu}\int_{0}^x (x_e-y)h(y)dy,$$
hence:
$$ (x_e-x)+\frac{1}{\mu}(x_e-x) \int_{0}^x h(y)dy =\frac{1}{\mu}\int_{0}^x (x_e-y)h(y)dy,$$
and so:
$$ (x_e-x) =\frac{1}{\mu}\int_{0}^x (x-y)h(y)dy,$$
which gives:
\begin{equation}\label{xstar00}
x_e =\frac{1}{\mu}\int_{0}^x (x-y)h(y)dy + x.
\end{equation}

Thus $x^*$ is a solution to~(\ref{xstar00}), and so: 
 \ $x^*=f^{-1}(x_e).$ \ It is easy to verify that $x^*$ is indeed a maximum point.


 Now, $f$ is strictly increasing since for all $s>t:$ \ $$f(s)-f(t)\geq\frac{1}{\mu}\int_{t}^s (s-y)h(y)dy + s-t\geq s-t>0,$$ and since $f(x)\geq x, \ \forall{x},$  then
 $\lim_{x\to\infty}f(x)=\infty.$ Since  \ $f(0)=0,$ then for every $x_e> 0,$ \ $x^*> 0.$
 In addition, because $h$ is locally integrable then by the fundamental theorem of calculus for Lebesgue integrals~\cite{Folland} $f$ is continuous.
It follows that for every $x_e\geq 0,$  $x^*$  is unique  and we have $\lim\limits_{x_e\to\infty}f^{-1}(x_e)=\infty.$ 

\hfill\quad$\Box$

In the sequel we will not refer to $f,f',$ hence  we summarize the above results in terms of $x_e, x^*$ as follows:
\begin{equation}\label{xstar1}
\lim_{x_e\to\infty}{x^*}=\infty.
\end{equation}
\begin{equation}\label{einc}
x^* \textsl{strictly increasing in} \  x_e.
\end{equation}
\begin{equation}\label{xstar0}
x_e =\frac{1}{\mu}\int_{0}^{x^*} (x^*-y)h(y)dy + x^*.
\end{equation}
Note that~(\ref{xstar0}) implies:
\begin{equation}\label{xstar2}
x_e \ge x^*.
\end{equation}

\

Note that the ``selfish'' optimal strategy $x_e,$  does not depend on the function $h.$ This is not surprising as the selfish customer
enters the empty system whenever the benefit exceeds the costs, and this does not depend on the arrival of other customers. In contrast,  $x^*,$ which is the optimal social strategy depends both on $x_e$ and on $h,$ as can be seen in~(\ref{xstar0}). It increases in $x_e,$ and decreases in $h,$ namely, if $h_1(x)\leq h_2(x), \forall{x\geq 0},$ then for any given  $x_e,$ the $x^*$ obtained in~(\ref{xstar0}) for $h_1$ is larger than the $x^*$ obtained in~(\ref{xstar0}) for $h_2.$ To understand this intuitively, note that a larger $x_e$ reflects higher utility of service, which makes $x^*$ larger as well.  Larger $h$ means that there are more potential customers (for any given threshold $x$), and therefore to gain a specific social benefit one needs a smaller $x^*.$

Define the price of anarchy $PoA(h,x_e)$ as: $$PoA(h,x_e)=\frac{S(x^*)}{S(x_e)}.$$
\

Since $x^*$ is the maximum point of $S(x)$ then \  $PoA(h,x_e)\geq 1.$

Now:
\begin{equation}\label{generalpoa}
PoA(h,x_e)=\frac{ \frac{c_t\left(\int_{0}^{x^*} (x_e-y)h(y)dy \right)}{1+\frac{1}{\mu} \int_{0}^{x^*} h(y)dy}}{ \frac{c_t\left(\int_{0}^{x_e} (x_e-y)h(y)dy \right)}{1+\frac{1}{\mu} \int_{0}^{x_e} h(y)dy}}=\left(
 \frac{\int_{0}^{x^*} (x_e-y)h(y)dy }{\int_{0}^{x_e} (x_e-y)h(y)dy}\right) \left(\frac{1+\frac{1}{\mu} \int_{0}^{x_e} h(y)dy}{1+\frac{1}{\mu} \int_{0}^{x^*} h(y)dy}\right).
 \end{equation}

\

Recall that $x_e$ is determined by the parameters of the model. In particular,
since $x_e=\frac{R\mu -c_w}{c_t\mu},$ then if $R\to\infty,$ or $c_t\to 0,$ we get $x_e\to\infty,$ and if   \  $R\mu\to c_w,$ or $c_t\to\infty,$ then $x_e\to 0.$
We wish to analyse the behaviour of the price of anarchy when $x_e\to 0,$ \ and when $x_e\to\infty.$ The first case is simple and is considered in the  following proposition. The rest of the paper analyses the behaviour of  $PoA(h, x_e)$  when $x_e\to\infty.$

\begin{Pro}\label{Tzero}
For all intensity functions $h,$ the price of anarchy $PoA(h,x_e)$ satisfies :

$$\lim\limits_{x_e\to 0}PoA(h, x_e)=1.$$
\end{Pro}

\proof

Note that the first factor appearing in the right-hand side of \eqref{generalpoa} is always smaller than $1$ (as $x_e\geq x^*$), and  the second factor converges to $1$ when $x_e \to 0,$ (since by  \eqref{xstar0}, $x^*\to 0$ when $x_e \to 0$). Since we always have \ $PoA(h,x_e)\geq 1,$ then
 $\lim\limits_{x_e\to 0}PoA(h, x_e)=1.$

\hfill\quad$\Box$


We turn now to analyse the limit of the price of anarchy when $x_e\to\infty.$
\

Recall that our process is defined such that for all $x\geq 0,$ customers with distances less than $x,$ arrive to the system according to a Poisson process with rate $\lambda (x)=\int\limits_0^x h(y)dy.$
We now show that there is a clear difference between cases in which $\lim\limits_{x\to \infty}\lambda (x)=\int\limits_0^{\infty} h(y)dy$ is finite and cases in which it is infinite.

\begin{Pro}\label{Pxstar}
\
\begin{itemize}
\item If $\int_0^{\infty} h(y) dy =\infty$ \ then \ $\lim\limits_{x_e\to \infty}\frac{x^*}{x_e}=0,$\notag\\
\item If $\int_0^{\infty} h(y) dy =a<\infty$ \ then \ $\frac{x^*}{x_e}\geq\frac{\mu}{a+\mu},$ for all $x_e.$\notag
\end{itemize}
\end{Pro}

\proof
By~(\ref{xstar0}):
\begin{equation}\label{EP}
\frac{x^*}{x_e}=\frac{x^*}{\frac{1}{\mu}\int_0^{x^*}(x^*-y)h(y)dy + x^*}=\frac{1}{\frac{1}{\mu}\int_0^{x^*}(1-\frac{y}{x^*})h(y)dy + 1}.
\end{equation}

Consider first the case in which $\int_0^{\infty} h(y) dy =\infty.$
Note that: $$\int_0^{x^*}(1-\frac{y}{x^*})h(y)dy\geq \int_0^{\frac{x^*}{2}}(1-\frac{y}{x^*})h(y)dy\geq\frac{1}{2}\int_0^{\frac{x^*}{2}} h(y) dy.$$

Recall that $\lim\limits_{x_e\to \infty}{x^*}=\infty,$ hence when $x_e\to\infty,$ the limit of the right hand side of this inequality  is infinite. Thus the limit of the left hand side is infinite as well, and so returning to~(\ref{EP}) we get  in this case: $\frac{x^*}{x_e}\to 0.$

\

If on the other hand,  $\int_0^{\infty} h(y) dy=a<\infty$ then:
$$\int_0^{x^*}(1-\frac{y}{x^*})h(y)dy\leq \int_0^{x^*} h(y) dy \leq a,$$

and so by~(\ref{EP}), for any given $x_e:$ \  $\frac{x^*}{x_e}\geq \frac{1}{\frac{a}{\mu}+1}=\frac{\mu}{a+\mu}.$

\hfill\quad$\Box$

The following theorem is interesting, as it suggests that according to the anarchy function, if the amount of customers arriving from far is ``small'', in the sense that the rate of arrival is bounded, then in the limit there is no difference between the social and equilibrium optimal benefits even though the corresponding optimal strategies $x^*$ and $x_e$ do not coincide.

\begin{Th}\label{Tfin}
For all intensity functions $h,$ if $\int_0^{\infty} h(y) dy =a<\infty,$ then:
$$\lim\limits_{x_e\to \infty}PoA(h, x_e)=1.$$
\end{Th}

\proof
From (\ref{generalpoa}) we have:
\begin{equation}\label{E77}
PoA(h, x_e)= \left(
 \frac{\int_{0}^{x^*} (x_e-y)h(y)dy }{\int_{0}^{x_e} (x_e-y)h(y)dy}\right) \left(\frac{1+\frac{1}{\mu} \int_{0}^{x_e} h(y)dy}{1+\frac{1}{\mu} \int_{0}^{x^*} h(y)dy}\right).
\end{equation}

As explained earlier $x^*\to\infty$ as $x_e\to\infty,$  hence the second factor in the right hand side of the equation above converges to $1.$

Now, the first factor equals:
\begin{equation}\label{E88}
 \frac{\int_{0}^{x^*} (1-\frac{y}{x_e})h(y)dy }{\int_{0}^{x_e} (1-\frac{y}{x_e})h(y)dy}.
 \end{equation}
Let \ $I_{[0,x^*]}$ be the indicator function
 of $[0,x^*].$
Then for every fixed $y,$
$$I_{[0,x^*]}(1-\frac{y}{x_e})h(y)$$
is increasing in $x_e$ and converges to $h(y)$ when $x_e\to\infty.$ (Recall that $x^*$ is an increasing function of $x_e$).

Hence, by the monotone convergence theorem:
$$\lim\limits_{x_e\to \infty}\int_{0}^{x^*} (1-\frac{y}{x_e})h(y)dy=\lim\limits_{x_e\to \infty}\int_{0}^{\infty}I_{[0,x^*]}(1-\frac{y}{x_e})h(y) dy=\int_{0}^{\infty}h(y) dy=a.$$
Similarly for the denominator in~(\ref{E88}) (with $x_e$ instead of $x^*$), and so the first factor in~(\ref{E77}) converges to $1$ as well.
\hfill\quad$\Box$

\begin{Co}
If $h(y)=\beta y^{\alpha},$ \ $\alpha < -1, \ \beta>0$ then:
$$\lim\limits_{x_e\to \infty}PoA(h,x_e)=1.$$
\end{Co}

\

Theorem~\ref{Tfin} completes our analysis of the case in which $\int_0^{\infty} h(y) dy $ \ is finite.

\
From now on, we always assume that $\int_0^{\infty} h(y) dy =\infty.$


\




In the following theorem and corollary we present two formulas to calculate the limit of PoA (if exists).

\begin{Th}\label{T99}
If $\int_0^{\infty} h(y) dy =\infty$ then:
$$\lim\limits_{x_e\to \infty}PoA(h, x_e)=\lim_{x_e\to\infty}\frac{ \int_{0}^{x_e} h(y)dy}{\int_{0}^{x_e} (1-\frac{y}{x_e})h(y)dy},$$
by which we mean that the first limit exists \emph{iff} the second limit exists and in that case they are equal.
\end{Th}

\proof
From (\ref{generalpoa}) we have:
\begin{equation}\label{EIM}
PoA(h, x_e)=
\left(
 \frac{\int_{0}^{x^*} (1-\frac{y}{x_e})h(y)dy }{1+\frac{1}{\mu} \int_{0}^{x^*} h(y)dy}\right) \left(\frac{1+\frac{1}{\mu} \int_{0}^{x_e} h(y)dy}{\int_{0}^{x_e} (1-\frac{y}{x_e})h(y)dy}\right).
 \end{equation}

Now, since for all $0\leq y\leq x^*:$
$$1-\frac{x^*}{x_e}\leq 1-\frac{y}{x_e}\leq 1,$$ then:
$$(1-\frac{x^*}{x_e})\int_0^{x^*} h(y) dy\leq \int_0^{x^*}(1-\frac{y}{x_e}) h(y) dy\leq \int_0^{x^*} h(y) dy.$$

So:
\begin{equation}\label{EIM1}
\frac{(1-\frac{x^*}{x_e})\int_0^{x^*} h(y) dy}{{1+\frac{1}{\mu} \int_{0}^{x^*} h(y)dy}}\leq \frac{\int_0^{x^*}(1-\frac{y}{x_e}) h(y) dy}{{1+\frac{1}{\mu} \int_{0}^{x^*} h(y)dy}}\leq\frac{\int_0^{x^*} h(y) dy}{{1+\frac{1}{\mu} \int_{0}^{x^*} h(y)dy}}.
\end{equation}
Now,  since $\lim\limits_{x_e\to \infty}\int_0^{x^*} h(y) dy =\infty$ (because $\lim\limits_{x_e\to \infty}x^*=\infty$) then by Proposition~\ref{Pxstar} $\lim\limits_{x_e\to \infty}\frac{x^*}{x_e}=0.$ It follows that the expressions appearing on the left- and on the right-hand sides of~(\ref{EIM1}) both converge to $\mu.$ Thus also the expression in the middle converges to $\mu,$ namely:
\begin{equation}\label{EIM3}
\lim_{x_e\to\infty}\frac{\int_0^{x^*}(1-\frac{y}{x_e}) h(y) dy}{{1+\frac{1}{\mu} \int_{0}^{x^*} h(y)dy}}=\mu.
\end{equation}

Returning to Equation~(\ref{EIM}) we get:
\begin{equation}\label{EIM4}
\lim_{x_e\to\infty}PoA(h, x_e)=\lim_{x_e\to\infty}\frac{\mu \left( 1+\frac{1}{\mu} \int_{0}^{x_e} h(y)dy\right)}{\int_{0}^{x_e} (1-\frac{y}{x_e})h(y)dy}.
\end{equation}

Now: \begin{equation}\label{E66}
\lim\limits_{x_e\to \infty}\int_{0}^{x_e} (1-\frac{y}{x_e})h(y)dy\geq \lim\limits_{x_e\to \infty}\int_{0}^{\frac{x_e}{2}} (1-\frac{y}{x_e})h(y)dy\geq \lim\limits_{x_e\to \infty}\frac{1}{2}\int_{0}^{\frac{x_e}{2}}h(y)dy=\infty.
\end{equation}
By this we get from~(\ref{EIM4}):

\begin{equation}\label{EIM5}
\lim_{x_e\to\infty}PoA(h, x_e)=\lim_{x_e\to\infty}\frac{ \int_{0}^{x_e} h(y)dy}{\int_{0}^{x_e} (1-\frac{y}{x_e})h(y)dy},
\end{equation}
by which we mean that the first limit exists \emph{iff} the second limit exists and in that case they are equal.

\hfill\quad$\Box$

The following corollary gives another presentation for \ $\lim_{x_e\to\infty}PoA(h, x_e).$
\begin{Co}\label{PANO}
If $\int_0^{\infty} h(y) dy =\infty$ then:
$$\lim\limits_{x_e\to \infty}PoA(h, x_e)=\lim\limits_{x_e\to \infty}\frac{\int_{0}^{1}h(x_et)dt}{\int_{0}^{1}(1-t)h(x_et)dt},$$
by which we mean that the first limit exists \emph{iff} the second one exists and in that case they are equal.
\end{Co}
\proof
Substituting $y=x_et,\ dy=x_edt$ in Theorem~\ref{T99} we get the required statement.
\hfill\quad$\Box$

\

Note that in Theorem~\ref{T99} and Corollary~\ref{PANO} the limit of the anarchy function when $x_e\to\infty$ is expressed without $x^*.$
In what follows, we always use the formula appearing in Theorem~\ref{T99}, except for Example~\ref{EXBOUND}, in which we use the formula introduced in Corollary~\ref{PANO}.

\begin{Pro}\label{PINV}
For all intensity functions $h,$ and for all $b>0:$
$$\lim_{x_e\to\infty}PoA(bh, x_e)=\lim_{x_e\to\infty}PoA(h, x_e),$$
by which we mean that the first limit exists \emph{iff} the second limit exists and in that case they are equal.
\end{Pro}

\proof
If  $\int_0^{\infty} h(y) dy =\infty$ then
By Corollary~\ref{PANO}:
$$\lim\limits_{x_e\to \infty}PoA(bh, x_e)=\lim\limits_{x_e\to \infty}\frac{\int_{0}^{1}bh(x_et)dt}{\int_{0}^{1}(1-t)bh(x_et)dt}=\lim\limits_{x_e\to \infty}\frac{\int_{0}^{1}h(x_et)dt}{\int_{0}^{1}(1-t)h(x_et)dt}=\lim\limits_{x_e\to \infty}PoA(h, x_e).$$

If however \ $\int_0^{\infty} h(y) dy$ is finite then also: $\int_0^{\infty} bh(y) dy$ is finite hence by Theorem~\ref{Tfin}:
$$\lim_{x_e\to\infty}PoA(bh, x_e)=\lim_{x_e\to\infty}PoA(h, x_e)=1.$$
\hfill\quad$\Box$

\begin{Co}\label{Con}
If $h(y)=\beta y^{\alpha},$ \ $\alpha > -1, \ \beta>0$ then:
$$\lim\limits_{x_e\to \infty}PoA(h,x_e)=\alpha+2.$$
In particular:
If $h$ is constant, namely: $h\equiv c,$ \ $c>0,$ then:
$$\lim\limits_{x_e\to \infty}PoA(h,x_e)=2.$$
\end{Co}

\begin{Ex}\label{EX1}\hfill\quad
\normalfont

Consider first the one-dimensional city. Assume that the arrival rate is uniform, s.t., $h(y)= \lambda>0,$ for all $y\geq 0$. By Corollary~\ref{Con}, the price of anarchy in this case converges to $2,$ as $x_e\to\infty.$ In the case of two-dimensional city, with a uniform arrival:  $h(y)=2\lambda\pi y,$ $\forall y\geq 0,$ by Corollary~\ref{Con}, the price of anarchy converges to $3.$
Intuitively, the difference between the anarchy in the uniform $1$- and $2$-dimensional cases, is due to the fact that on the two-dimensional city there exist relatively more customers with high distances that enter the system. We argue that because of this difference, the anarchy is larger in the case of $2$-dimensional city than in the $1$-dimensional city. This intuitive argument can be generalised for all  uniform $n$-dimensional city, to explain why we get $\lim\limits_{x_e\to \infty}PoA(x^n,x_e)=n+2,$ which increases in $n.$
Thus the higher the dimension, the higher the price of anarchy.
Later on we generalize this argument in Theorem~\ref{TTWO} to  $h_1, h_2,$ which satisfy  $h_1(x)\leq h_2(x)$ from some point on. We prove that under certain conditions: $\lim\limits_{x_e\to \infty}PoA(h_1, x_e)\leq\lim\limits_{x_e\to \infty}PoA(h_2, x_e)$ (see Theorem~\ref{TTWO}).
\end{Ex}

\begin{Th}\label{TDIV}
Given two intensity functions $h_1, h_2$ s.t. \  $\int_0^{\infty} h_1(y) dy =\int_0^{\infty} h_2(y) dy=\infty,$
if $$\lim\limits_{x\to \infty}\frac{h_1(x)}{h_2(x)}=c>0,$$ then:
$$\lim\limits_{x_e\to \infty}PoA(h_1, x_e)=\lim\limits_{x_e\to \infty}PoA(h_2, x_e),$$
namely, if one of the limits exists then so does the other, and in this case they are equal.
\end{Th}

\proof
Because of Proposition~\ref{PINV} we can replace $h_2$ by $ch_2$ and have: $\lim\limits_{x\to \infty}\frac{h_1(x)}{h_2(x)}=1.$

Define $u=h_1-h_2.$ Then:
$$\lim\limits_{x\to \infty}\frac{u(x)}{h_2(x)}=\lim\limits_{x\to \infty}\frac{h_1(x)-h_2(x)}{h_2(x)}=\lim\limits_{x\to \infty}\frac{h_1(x)}{h_2(x)}-1=0.$$


Now:
$$\frac{ \int_{0}^{x} h_1(y)dy}{\int_{0}^{x} (1-\frac{y}{x})h_1(y)dy}=\frac{ \int_{0}^{x}\left( h_2(y)+u(y)\right)dy}{\int_{0}^{x} (1-\frac{y}{x})\left( h_2(y)+u(y)\right) dy}=\frac{ \int_{0}^{x} h_2(y)dy+\int_{0}^{x}u(y)dy}{\int_{0}^{x} (1-\frac{y}{x}) h_2(y)+\int_{0}^{x}(1-\frac{y}{x})u(y) dy}$$

dividing the numerator and the denominator by $\int_0^xh_2(y)dy$ to get:
\begin{equation}\label{E55}
\frac{1+\frac{\int_0^xu(y)dy}{\int_0^xh_2(y)dy}}{\frac{\int_{0}^{x} (1-\frac{y}{x}) h_2(y)dy}{\int_0^xh_2(y)dy}+\frac{\int_{0}^{x}(1-\frac{y}{x})u(y) dy}{\int_0^xh_2(y)dy}}.
\end{equation}

Now, $\int_0^x h_2(y)dy\to\infty$ and so we can use L'Hopital's rule~\cite{LHopital} and get that: $\frac{\int_0^xu(y)dy}{\int_0^xh_2(y)dy}$ appearing in the numerator of~(\ref{E55}) converges to: $\lim\limits_{x\to \infty}\frac{u(x)}{h_2(x)}=0,$ so the numerator converges to $1.$

In addition: $\frac{\int_{0}^{x}(1-\frac{y}{x})u(y) dy}{\int_0^xh_2(y)dy}$ appearing in the denominator of~(\ref{E55}) satisfies: $$0\leq \frac{\int_{0}^{x}(1-\frac{y}{x})u(y) dy}{\int_0^xh_2(y)dy}\leq \frac{\int_0^xu(y)dy}{\int_0^xh_2(y)dy}\to 0,$$ and so it too converges to $0.$   Note that the limit of \
$\frac{\int_{0}^{x} (1-\frac{y}{x}) h_2(y)dy}{\int_0^xh_2(y)dy}$ appearing in the denominator of~(\ref{E55}) is  $\frac{1}{\lim\limits_{x_e\to \infty}PoA(h_2, x_e)}.$ Taken together, we have:
$$\lim\limits_{x_e\to \infty}PoA(h_1, x_e)=\lim\limits_{x_e\to \infty}PoA(h_2, x_e).$$

\hfill\quad$\Box$

\begin{Pro}\label{PFIX}
If $\lim\limits_{x\to \infty}h(x)=c>0,$  then: $$\lim\limits_{x_e\to \infty}PoA(h, x_e)=2.$$
\end{Pro}

\proof
First, $\lim\limits_{x\to \infty}h(x)=c,$  implies $\int_0^{\infty} h(y) dy=\infty.$
By Corollary~\ref{Con} $\lim\limits_{x_e\to \infty}PoA(\bar{1}, x_e)=2,$ where $\bar{1}$ is the constant function $1.$  Substituting $h_2(x)=\bar{1},$ in Theorem~\ref{TDIV} proves the proposition.

\hfill\quad$\Box$

\begin{Co}\label{CCIN}
If two intensity functions $h_1, h_2,$ coincide from some point on, then:
$$\lim\limits_{x_e\to \infty}PoA(h_1, x_e)=\lim\limits_{x_e\to \infty}PoA(h_2, x_e),$$
by which we mean that the first limit exists \emph{iff} the second limit exists and in that case they are equal.
\end{Co}

\proof
The proof follows immediately from Theorem~\ref{TDIV}.

\begin{Th}\label{p33}
If $\int_0^{\infty} h(y) dy =\infty$  and if  $\lim\limits_{x_e\to \infty}PoA(h, x_e)$ exists then:
\begin{enumerate}
\item If $h$ increases monotonically from some point on then:  $\lim\limits_{x_e\to \infty}PoA(h, x_e)\geq 2.$\label{L2}\item If $h$ decreases monotonically from some point on then:  $\lim\limits_{x_e\to \infty}PoA(h, x_e)\leq 2.$\label{L3}
\end{enumerate}
\end{Th}

\proof

Note first that it is sufficient to prove the theorem for monotonic functions, and then by Corollary~\ref{CCIN} the theorem is proved also for functions that are monotonic only from some point on.

Assume that $h$ is monotonically increasing then we will show that for all $x>0:$
$$\frac{\int_{0}^{x} h(y)dy}{\int_{0}^{x} (1-\frac{y}{x})h(y)dy}\geq 2,$$ and then by Theorem~\ref{T99} this will prove~\ref{L2}.

For each $x>0$ define: $$c_{x}=\int_0^x h(y) dy.$$
Note that: $h(0)\leq \frac{c_x}{x}$ and also $h(x)\geq \frac{c_x}{x}.$ This is true since if $h(0)>\frac{c_x}{x},$ then since $h$ is increasing then $c_x=\int_0^x h(y) dy>\int_0^x\frac{c_x}{x}dy=c_x,$ which is a contradiction.  Similarly, if $h(x)< \frac{c_x}{x}$ then $c_x=\int_0^x h(y) dy<\int_0^x\frac{c_x}{x}=c_x,$ a contradiction. Since $h$ is increasing this implies that there exists $0\leq \theta\leq x,$  s.t:
\begin{itemize}
\item $h(y)\leq\frac{c_x}{x},$ \ \ \ $\forall{y\leq\theta}$
\item $h(y)\geq\frac{c_x}{x},$ \ \ \ $\forall{y\geq\theta}.$
\end{itemize}

Now:
$$\left( 1-\frac{\theta}{x}\right)\int_0^x\left( h(y)-\frac{c_x}{x}\right)dy=\left( 1-\frac{\theta}{x}\right)\left[ \int_0^x h(y)dy -\frac{1}{x}\int_0^xc_x dy \right]=\left( 1-\frac{\theta}{x}\right)\left[ c_x-c_x\right]=0.$$
Hence:
\begin{equation}\label{E800}
\int_{0}^{x} \left( 1-\frac{y}{x}\right)\left(h(y)-\frac{c_x}{x}\right) dy=\int_{0}^{x} \left[\left(1-\frac{y}{x}\right)- \left(1-\frac{\theta}{x}\right) \right]\left(h(y)-\frac{c_x}{x}\right) dy=
\int_{0}^{x}\frac{\theta-y}{x}\left(h(y)-\frac{c_x}{x}\right) dy=$$
$$\int_{0}^{\theta}\frac{\theta-y}{x}\left(h(y)-\frac{c_x}{x}\right) dy+\int_{\theta}^{x}\frac{\theta-y}{x}\left(h(y)-\frac{c_x}{x}\right) dy.
\end{equation}

The integrand in the first integral is negative, being a product of a function which is positive and a function that is negative, in the given domain. The integrand of the second integral is similarly negative.
 Thus:
\begin{equation}\label{E900}
\int_{0}^{x}\left( 1-\frac{y}{x} \right)\left( h(y)-\frac{c_x}{x}\right) dy\leq 0,
\end{equation}
and so:
\begin{equation}\label{E700}
\int_{0}^{x}\left( 1-\frac{y}{x} \right)h(y) dy\leq \int_{0}^{x}\left( 1-\frac{y}{x} \right)\frac{c_x}{x} dy=\frac{c_x}{x}\left( x-\frac{x}{2} \right)=\frac{c_x}{2}=\frac{\int_0^x h(y) dy}{2},
\end{equation}
hence:
 $$\frac{\int_{0}^{x} h(y)dy}{\int_{0}^{x} (1-\frac{y}{x})h(y)dy}\geq 2.$$

This proves~\ref{L2}.
The proof of~\ref{L3} is  similar.
\hfill\quad$\Box$

Note that the limit of $PoA(h, x_e)$ does not always exist.
In particular, this implies that PoA is not always monotone in $x_e.$
The following two examples show two cases in which $\lim\limits_{x_e\to \infty}PoA(h,x_e)$ does not exist. In the first example $h$ is monotonic and in the second example, $h$ is bounded between two  strictly positive constants.

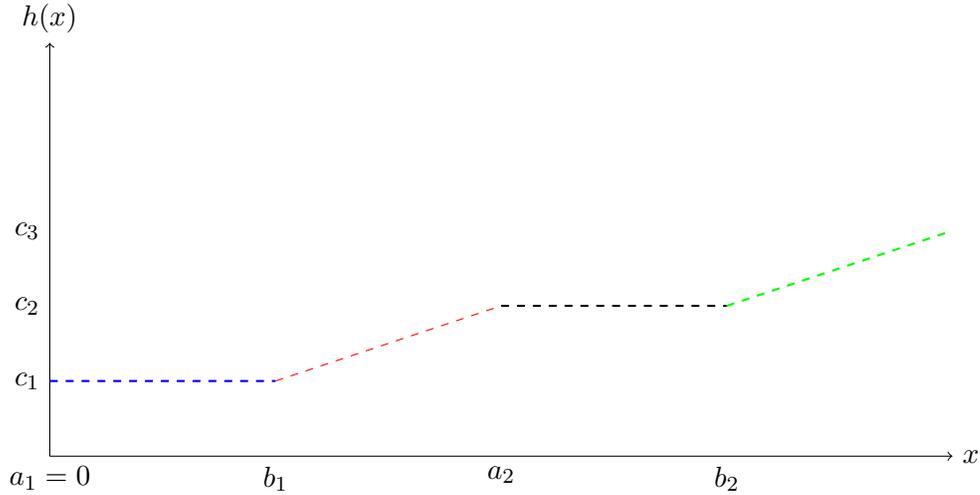
\begin{figure}[H]
\centering
\begin{tikzpicture}[xscale=3,yscale=1]
  \def\xmin{0}
  \def\xmax{4}
  \def\ymin{0}
  \def\ymax{5.5}
    \draw[->] (\xmin,\ymin) -- (\xmax,\ymin) node[right] {$x$} ;
    \draw[->] (\xmin,\ymin) -- (\xmin,\ymax) node[above] {$h$($x$)} ;
      \node at (0,0) [below] {$a_1=0$};
      \node at (1,0) [below] {$b_1$};
      \node at (2,0) [below] {$a_2$};
      \node at (3,0) [below] {$b_2$};
      \node at (0,1) [left] {$c_1$};
      \node at (0,2) [left] {$c_2$};
      \node at (0,3) [left] {$c_3$};
    \draw[blue,dashed, thick,domain=0:1]  plot (\x, {\x^0});
    \draw[red, dashed,domain=1:2]  plot (\x, {\x});
    \draw[black,dashed, thick,domain=2:3]  plot (\x, {2});
    \draw[green,dashed, thick,domain=3:4]  plot (\x, {\x-1});

        \end{tikzpicture}
\caption{An illustration of $h$ for Example~\ref{EX100}.}\label{F10}
\end{figure}

\begin{Ex}\label{EX100}
\normalfont

\

We define $h$ as a piecewise linear function (see Figure~\ref{F10}) as follows:
 Divide the domain $R^+$ into intervals $[a_i, b_i],$ and $[b_i, a_{i+1}], \ i=1,2,\dots $
Define: $$h(x)=c_i,  \  \ \ \  \forall{a_i\leq x\leq b_i},$$ and \ $$h(x)=(c_i-b_i)+x,  \ \ \ \ \forall{b_i\leq x\leq a_{i+1}}.$$ Define $c_{i+1}=c_i-b_i+a_{i+1}.$  It is easy to verify that $h$ is continuous.
\

We now need to specify $a_i, b_i.$
Recall that when $h(x)= c_i$ for all  $x$ large enough, then $\lim\limits_{x_e\to \infty}PoA(h,x_e)=2,$ (Proposition~\ref{PFIX}) and when $h$ is linear from some point on, then $\lim\limits_{x_e\to \infty}PoA(h,x_e)=3,$ (Example~\ref{EX1} and Corollary~\ref{CCIN}).
\
Define $a_1=0.$ On $[0,b_1]:$ \ $h(x)=c_1, \ \forall{x}.$ Define: $h_1(x)$ as the function that equals $c_1$ for all $x\geq 0,$ namely $h_1=\bar{c_1}.$

By Corollary~\ref{Con} $\lim\limits_{x\to \infty}PoA(\bar{c_1}, x)=2,$ hence there exists $b_1$ s.t for all $x\geq b_1,$  $PoA(\bar{c_1}, x)\leq 2.3.$
Now, in $[b_1, a_2]$ \ $h$ is linear. Define: $h_2(x)$ as the function that equals $c_1$ on $[0, b_1]$ and from that point on equals $(c_2-b_2)+x.$   By Corollary~\ref{CCIN} $\lim\limits_{x\to \infty}PoA(h_2, x)=3,$ hence there exists $a_2\geq b_1$ s.t for  $x\geq a_2,$  $PoA(h_2, x)\geq 2.7,$  and so on we continue defining $a_i, b_i$ in the same way such that $h= h_i$ for all $x$ between 0 and the end of the $i$th interval. Thus we get that $PoA(h, x)\leq 2.3$ and then $PoA(h, x)\geq 2.7$ alternately, and so
 $PoA(h,x_e)$ does not have a limit.
\end{Ex}

\

\begin{Ex}\label{EXBOUND}
\normalfont

\

Divide the domain $R^+$ into intervals: $[0, 1), [1, 10),[10, 100),[100, 1000),\dots$

Define $h$ as the function that attains the constant functions $\bar{1}$ and $\bar{2}$ alternately on these intervals starting with $\bar{1}$ at the first interval $[0, 1).$

One can verify easily that for $x=2, 200, 20000\dots :$
$$\frac{\int_{0}^{1}h(xt)dt}{\int_{0}^{1}(1-t)h(xt)dt}\leq \frac{\frac{3}{2}}{\frac{7}{8}-\frac{1}{20}}=\frac{20}{11},$$

and for $x=20, 2000, 200000\dots :$
$$\frac{\int_{0}^{1}h(xt)dt}{\int_{0}^{1}(1-t)h(xt)dt}\geq \frac{\frac{3}{2}}{\frac{5}{8}+\frac{1}{20}}=\frac{20}{9}.$$

Thus according to Corollary~\ref{PANO} $PoA(h,x_e)$ does not have a limit. Note that $h$ defined in this example can be smoothed out to give a continuous function with similar properties.
\end{Ex}

\

The following theorem shows that in many cases the limit of the price of anarchy does exist and has a simple presentation.

\begin{Th}\label{TEX}
If $\int_0^{\infty} h(y) dy =\infty$ and $\lim\limits_{x\to \infty}\frac{xh'(x)}{h(x)}$ \ exists then $\lim\limits_{x_e\to \infty}PoA(h, x_e)$ exists and:
$$\lim\limits_{x_e\to \infty}PoA(h, x_e)=2+\lim\limits_{x\to \infty}\frac{xh'(x)}{h(x)}.$$
\end{Th}

\proof
From~(\ref{E66}) we have that if $\int_0^{\infty} h(y) dy =\infty,$ then also: $\lim\limits_{x_e\to \infty}\int_{0}^{x_e} (1-\frac{y}{x_e})h(y)dy=\infty,$ and so we can use L'Hopital's rule to get:
$$\lim_{x\to\infty}\frac{\int_{0}^{x} h(y)dy}{\int_{0}^{x} (1-\frac{y}{x})h(y)dy}=\lim_{x\to\infty}\frac{h(x)}{\int_{0}^{x} \frac{y}{x^2}h(y)dy}=\lim_{x\to\infty}\frac{x^2h(x)}{\int_{0}^{x}yh(y)dy}.$$

(the contribution of the derivative with respect to the upper limit of the integral appearing in the denominator on the left-hand side of the equation above is $0,$ since $\left( 1-\frac{y}{x}\right)h(y)$ is $0$ for $y=x$).

Using L'Hopital's rule again we get that this equals:
$$\lim\limits_{x\to \infty}=\frac{2xh(x)+x^2h'(x)}{xh(x)}=2+\lim\limits_{x\to \infty}\frac{xh'(x)}{h(x)}.$$

\hfill\quad$\Box$
\

Note that the number $2$ plays a central role in our work (see Theorems~\ref{p33} and~\ref{TEX}, Corrolary~\ref{Con} and Proposition~\ref{PFIX}). This is interesting since as mentioned earlier, Gilboa-Freedman, Hassin and Kerner~\cite{GHK} showed that if $\lambda=\mu$ in Naor's model, then the price of anarchy equals $2$ as well.

\begin{Exs}\label{Efinal}\hfill\quad
\normalfont
\begin{itemize}
\item   $h(x)=e^x.$

Then $h'(x)=e^x$ as well, and so:
$$\lim\limits_{x\to \infty}\frac{xh'(x)}{h(x)}=\lim\limits_{x\to \infty} x=\infty,$$
hence by Theorem~\ref{TEX}:
$$\lim\limits_{x_e\to \infty}PoA(e^x, x_e)=2+\infty=\infty.$$
\item $h(x)=\ln{(x+1)}.$

Then $h'(x)=\frac{1}{x+1},$ and so:
$$\lim\limits_{x\to \infty}\frac{xh'(x)}{h(x)}=\lim\limits_{x\to \infty} \left(\frac{x}{x+1}\right)\left(\frac{1}{\ln{(x+1)}}\right)=0,$$
hence by Theorem~\ref{TEX}:
$$\lim\limits_{x_e\to \infty}PoA(\ln{(x+1)}, x_e)=2+0=2.$$
\item $h(x)=\frac{1}{x+1}.$

Then $h'(x)=-\frac{1}{(x+1)^2},$ and so:
$$\lim\limits_{x\to \infty}\frac{xh'(x)}{h(x)}=\lim\limits_{x\to \infty}-\left(\frac{x}{x+1}\right)=-1,$$
hence by to Theorem~\ref{TEX}:
$$\lim\limits_{x_e\to \infty}PoA\left(\frac{1}{x+1}, x_e\right)=2+(-1)=1.$$

The following  example shows a case in which we cannot use Theorem~\ref{TEX} since  $\lim\limits_{x\to \infty}\frac{xh'(x)}{h(x)}$ does not exist. In such a case $\lim\limits_{x_e\to \infty}PoA(h, x_e)$ should be computed directly from Theorem~\ref{T99} or Corollary~\ref{PANO}.

\item  $h(x)=2+\sin {x}.$

Then $h'(x)=\cos {x},$ and so:
$$\frac{xh'(x)}{h(x)}=\frac{x\cos {x}}{2+\sin {x}}.$$
To see that this expression does not have a limit (and therefore we cannot use Theorem~\ref{TEX}), denote: $u(x)=\frac{x\cos {x}}{2+\sin {x}}.$ Then: $u(2k\pi )=k\pi,$ \ whereas $u((2k+1)\pi )= -\left( k+\frac{1}{2}\right)\pi ,$ and so $\frac{xh'(x)}{h(x)}$ does not have a limit.
\
However according to Theorem~\ref{T99}:

$$\lim\limits_{x_e\to \infty}PoA(2+\sin {x}, x_e)=\lim_{x_e\to\infty}\frac{ \int_{0}^{x_e} (2+\sin {y}) dy}{\int_{0}^{x_e} (1-\frac{y}{x_e})(2+\sin {y}) dy}=\lim\limits_{x_e\to \infty}\frac{2x_e-\cos {x_e}+1}{x_e+1-\frac{1}{x_e}\sin {x_e}}$$
$$=\lim\limits_{x_e\to \infty}\frac{2-\frac{1}{x_e}\cos {x_e}+\frac{1}{x_e}}{1+\frac{1}{x_e}-\frac{1}{x_e^2}\sin {x_e}}=2.$$
\end{itemize}
\end{Exs}

The following theorem generalizes our observation regarding $h_{n}(y)=y^n$ (see Example~\ref{EX1}), that the more customers  arriving from far, the higher the limit of the price of anarchy.

\begin{Th}\label{TTWO}
Given two intensity functions $h_1,h_2$ with infinite integrals  and for which the limits of $\frac{xh'_i(x)}{h_i(x)}$  exist ($i=1,2$) then if there exists $M>0,$ s.t for all $x\geq M:$ $h_1,h_2$ are monotonic and $h_1(x)\leq h_2(x),$ then:
$$\lim\limits_{x_e\to \infty}PoA(h_1, x_e)\leq\lim\limits_{x_e\to \infty}PoA(h_2, x_e).$$
\end{Th}

\proof
If $h_1$ is decreasing  and  $h_2$ is increasing then by Theorem~\ref{p33} \ $$\lim\limits_{x_e\to \infty}PoA(h_1, x_e)\leq 2\leq \lim\limits_{x_e\to \infty}PoA(h_2, x_e).$$
If $h_1$ is increasing  and  $h_2$ is decreasing then both must converge to a positive constant and so by Proposition~\ref{PFIX} \ $\lim\limits_{x_e\to \infty}PoA(h_1, x_e)= 2=\lim\limits_{x_e\to \infty}PoA(h_2, x_e).$

Hence we need to prove the theorem for the case that the functions are either both increasing monotonically or both decreasing monotonically. If both are increasing and $h_2$ converges to a positive constant then so does $h_1$ and in that case by Proposition~\ref{PFIX} $\lim\limits_{x_e\to \infty}PoA(h_i, x_e)=2,$ for both $i=1,2.$ If only $h_1$ converges to a positive constant then since $h_2$ is increasing then according to  Proposition~\ref{PFIX} and
Theorem~\ref{p33}: $\lim\limits_{x_e\to \infty}PoA(h_1, x_e)= 2\leq \lim\limits_{x_e\to \infty}PoA(h_2, x_e).$ Hence if both are increasing we only need to prove the theorem for the case in which both functions are increasing monotonically to infinity. Similarly, if both are decreasing we need to prove the theorem only for the case in which they both converge to zero.

Assume that $h_1, h_2$ both decrease monotonically to $0.$ Then $h'_1, h'_2\leq 0,$ and so there exist $a,b$ s.t:
$$\lim\limits_{x\to \infty}\frac{xh'_1(x)}{h_1(x)}=a\leq 0,\ \ \emph{and} \ \  \lim\limits_{x\to \infty}\frac{xh'_2(x)}{h_2(x)}=b\leq 0.$$

Assume  contrary to our statement, that $\lim\limits_{x_e\to \infty}PoA(h_1, x_e)>\lim\limits_{x_e\to \infty}PoA(h_2, x_e).$ Then by Theorem~\ref{TEX} \ $b<a\leq 0,$ and so: \ $0\leq \frac{a}{b}<1,$ (since $b$ is negative).

Hence:
$$\lim\limits_{x\to \infty}\frac{\frac{h'_1(x)}{h_1(x)}}{\frac{h'_2(x)}{h_2(x)}}=\lim\limits_{x\to \infty}\frac{\frac{xh'_1(x)}{h_1(x)}}{\frac{xh'_2(x)}{h_2(x)}}=\frac{a}{b}<1.$$

In other words:
\begin{equation}\label{ELN}
\lim\limits_{x\to \infty}\frac{(\ln{h_1})'}{(\ln{h_2})'}=\frac{a}{b}<1.
\end{equation}

Now, since $h_2\to 0,$ then $\ln{h_2}\to -\infty,$ so we can use L'Hopital's rule in~(\ref{ELN}) to get:
$$\lim\limits_{x\to \infty}\frac{\ln{h_1}}{\ln{h_2}}=\frac{a}{b}<1.$$
Thus there exists $A>0,$ s.t. for all $x\geq A:$
$\frac{\ln {h_1}}{\ln{h_2}}<1.$ \ Since \  $\ln{h_2(x)}<0$ for $x$ large enough, then we get: $\ln{h_1(x)}>\ln{h_2(x)},$ implying: $h_1(x)>h_2(x),$ for all $x\geq A,$ contradicting the assumption of the theorem, that there exists $M>0,$ s.t for all $x\geq M:$
$h_1(x)\leq h_2(x)$.

The proof for the case in which $h_1$ and $h_2$ are both \emph{increasing} monotonically to infinity is  similar.

\hfill\quad$\Box$


\begin{Ex}\label{E900}
Return to $h(x)=e^x$ in Example~\ref{Efinal} and $h_n(x)=x^n$ discussed earlier in Example~\ref{EX1}.
Since for all $n$ \ $x^n<e^x$ from some point on, then according to Theorem~\ref{TTWO}:
$$\lim\limits_{x_e\to \infty}PoA(x^n, x_e)\leq\lim\limits_{x_e\to \infty}PoA(e^x, x_e), \ \forall{n}.$$
Indeed, as shown earlier: \ $\lim\limits_{x_e\to \infty}PoA(x^n, x_e)=n+2,$ and \ $\lim\limits_{x_e\to \infty}PoA(e^x, x_e)=\infty.$
\end{Ex}

\

\section{A system with a queue}

In this section, we analyse a system with a queue, with a uniform intensity function, s.t.,  $h(y)\equiv \lambda,  \forall{y>0.}$ The optimal strategy of a
customer located near the server (namely $x=0$) is to enter the system \textit{iff} the
queue length $n$ satisfies: $R \ge \frac{c_w(n+1)}{\mu}$, and so:
$n+1\le \frac{R\mu}{c_w}$, equivalently $n< n_e=\lfloor
\frac{R\mu}{c_w} \rfloor$ like in \cite{N}. A customer located at a distance $x$ from the server and
observing a queue length $i$, $0\le i\le n_e-1$ joins the system  \textit{iff}
$R-\frac{c_w(i+1)}{\mu}-c_tx\ge 0$.
Define $x_i^e$ as the Nash equilibrium threshold when the queue length is $i$. Hence the equilibrium strategy is characterized by a vector of thresholds
$(x_0^e, x_1^e,\ldots, x_{n_e-1}^e)$ where:
$$x_0^e=\frac{R\mu -
c_w}{c_t\mu}=\frac{\nu-1}{\kappa}, \quad \textit{for queue length} \quad 0,$$
$$x_1^e=\frac{R\mu -
2c_w}{c_t\mu}=\frac{\nu-2}{\kappa}, \quad \textit{for queue length} \quad 1,$$
$$\vdots
$$
$$x_{n_e-2}^e=\frac{R\mu -(n_e-1)c_w}{c_t\mu}=\frac{\nu-(n_e-1)}{\kappa}, \quad \textit{for queue length} \quad n_e-2,$$

\begin{equation}\label{E5}
x_{n_e-1}^e=\frac{R\mu -
n_ec_w}{c_t\mu}=\frac{\nu-n_e}{\kappa}, \quad \textit{for queue length} \quad n_e-1.
\end{equation}
Thus when a customer observes a queue length $i$, he joins the system \textit{if}
$x\le x_i^e$ and balks otherwise.

Recall that $x_e$ was defined earlier as the equilibrium strategy in the loss system and note that $x_0^e$ equals $x_e,$  as both equal $\frac{\nu -1}{\kappa}.$

Let $S(x_0, x_1,\ldots, x_{n_e-1})$ be the expected social benefit
per unit of time when the threshold strategy is $(x_0, x_1,\ldots,
x_{n_e-1})$.  We consider the case of one-dimensional uniform arrival. Namely, the arrival rates are $\lambda x_0, \lambda x_1,\ldots, \lambda x_{n_e-1}$
(depending on the queue length) and the probability vector of observing $n$
customers in the system $(n=0,1,...,n_e)$ is $(\pi_0, \pi_1,\ldots,\pi_{n_e})$.
Therefore, $(\pi_0, \pi_1,\ldots,\pi_{n_e})$ is a solution to the
balance equations system:
$$-\lambda x_0\pi_0+\mu\pi_1=0,$$
$$\lambda x_0\pi_0-(\lambda x_1+\mu)\pi_1+\mu\pi_2=0,$$
$$\lambda x_1\pi_1-(\lambda x_2+\mu)\pi_2+\mu\pi_3=0,$$
$$\vdots
$$
$$\lambda x_{n_e-1}\pi_{n_e-1}-\mu\pi_{n_e}=0.$$
The solution of the equations system is $\pi_i=\rho^i x_0\cdots
x_{i-1}\pi_0$, $1\le i\le n_e$ where $\rho=\frac{\lambda}{\mu},$ and:
 $$\pi_0=\frac{1}{1+\rho x_0+\cdots+\rho^{n_e}x_0\cdots x_{n_e-1}}.$$
The expected social benefit satisfies the following equations,
$$S(x_0, x_1,\ldots, x_{n_e-1})=\sum\limits_{n=1}^{n_e}\int_0^{x_{n-1}} \left(R -\frac{c_w n}{\mu}-c_ty\right)\lambda\pi_{n-1}(x_0, x_1,\ldots, x_{n_e-1})dy$$
$$=\sum\limits_{n=1}^{n_e} \left(R -\frac{c_w n}{\mu}-\frac{c_tx_{n-1}}{2}\right)x_{n-1}\lambda \pi_{n-1}(x_0, x_1,\ldots, x_{n_e-1})$$
$$=\sum\limits_{n=1}^{n_e} \left(x^e_{n-1}-\frac{x_{n-1}}{2}\right)c_tx_{n-1}\lambda\pi_{n-1}(x_0, x_1,\ldots, x_{n_e-1})$$
$$=\frac{\sum\limits_{n=1}^{n_e} \left(x^e_{n-1}-\frac{x_{n-1}}{2}\right)\lambda c_tx_{n-1} \rho^{n-1} x_0\cdots
x_{n-2}}{1+\rho x_0+\cdots+\rho^{n_e} x_0\cdots x_{n_e-1}}$$
$$=\frac{\sum\limits_{n=1}^{n_e} \left(x^e_{n-1}-\frac{x_{n-1}}{2}\right)\lambda c_t \rho^{n-1} x_0\cdots
x_{n-1}}{1+\sum\limits_{n=1}^{n_e} \rho^{n} x_0\cdots
x_{n-1}}.$$




In the loss system we proved that the price of anarchy is always bounded. The following theorem shows that this does not hold in the case of a system with a queue.
\begin{Th}\label{unbounded}
The price of anarchy, $PoA(\rho,x_0^e, x_1^e,\ldots, x_{n_e-1}^e),$ is unbounded.
\end{Th}

\proof
Recall that $x^*$ denotes the optimal solution for the loss system, and recall that: $x_0^e=x_e.$
Then: $$S(x_0^*, x_1^*,\ldots, x_{n_e-1}^*)\ge S(x^*, 0,0,\ldots, 0)=\frac{c_t x^*}{1+\rho x^*}\left(x_0^e-\frac{x^*}{2}\right).$$
The above is true since the right-hand side is the optimal social benefit $S(x^*)$ for the loss-system when $h(y)\equiv \lambda$ (see \eqref{socialbenefit}), and is not necessarily optimal for a system with a queue. Thus:
$$PoA(\rho, x_0^e, x_1^e,\ldots, x_{n_e-1}^e)=\frac{S(x_0^*, x_1^*,\ldots, x_{n_e-1}^*)}{S(x_0^e, x_1^e,\ldots, x_{n_e-1}^e)} \ge
\frac{\frac{c_t x^*}{1+\rho x^*}\left(x_e-\frac{x^*}{2}\right)}{\frac{\sum\limits_{k=1}^{n_e} \left(x^e_{k-1}-\frac{x_{k-1}^e}{2}\right)c_t
\rho^{k-1} x_0^e\cdots
x_{k-1}^e}{1+\rho x_0^e+\cdots+\rho^{n_e} x_0^e\cdots x_{n_e-1}^e }}$$

$$=\frac{\frac{ x^*}{1+\rho x^*}\left(x^e_0-\frac{x^*}{2}\right)}{\frac{\frac{1}{2}\sum\limits_{k=1}^{n_e} \left(x^e_{k-1}\right)
\rho^{k-1} x_0^e\cdots
x_{k-1}^e}{1+\rho x_0^e+\cdots+\rho^{n_e} x_0^e\cdots x_{n_e-1}^e }}=2 \cdot\frac{\frac{ x^*}{1+\rho x^*}\left(x^e_0-\frac{x^*}{2}\right)}{\frac{\sum\limits_{k=1}^{n_e}
\rho^{k-1} x_0^e\cdots
(x_{k-1}^e)^2}{1+\rho x_0^e+\cdots+\rho^{n_e} x_0^e\cdots x_{n_e-1}^e }}$$

$$=2 \cdot\frac{x^*}{(1+\rho x^*)} \cdot \frac{x^e_0-\frac{x^*}{2}}{x^e_0}\cdot \frac{x_0^e(1+\rho x_0^e+\cdots+\rho^{n_e} x_0^e\cdots
x_{n_e-1}^e)}{(x_0^e)^2+ \rho x_0^e(x_1^e)^2+\rho^2x_0^ex_1^e(x_2^e)^2+,\ldots,+ \rho^{n_e-1}x_0^e\cdots
x_{n_e-2}^e(x_{n_e-1}^e)^2}.$$

\begin{equation}\label{que}
=2 \cdot\Biggl(\frac{1}{\frac{1}{x^*}+\rho} \Biggr)\Biggl(1-\frac{1}{2}\Bigl(\frac{x^*}{x_e} \Bigr) \Biggr)\Biggl( \frac{1+\rho x_0^e+\cdots+\rho^{n_e} x_0^e\cdots
x_{n_e-1}^e}{x_0^e+ \rho (x_1^e)^2+\rho^2x_1^e(x_2^e)^2+,\ldots,+ \rho^{n_e-1}\cdots
x_{n_e-2}^e(x_{n_e-1}^e)^2}\Biggr)
\end{equation}


Denote $\sqrt{\frac{c_w}{\mu}}$ as  $s.$  Given $s,$ choose the parameters $R,$ and $c_t,$ such that:  $c_t=1,$  and $R=\frac{(2s-1)s^2}{s-1}.$
Substituting these choices in~(\ref{E5}), yields:

$$x_0^e=R-c_w/\mu=\frac{s}{s-1}\cdot s^2,$$
$$x_1^e=R-2c_w/\mu=\frac{1}{s-1}\cdot s^2.$$
which implies:
$$\frac{x_0^e}{x_1^e}=\frac{R-c_w/\mu}{R-2c_w/\mu}=s.$$

Recall that $n_e=\lfloor \frac{R\mu}{c_w}\rfloor=\lfloor \frac{2s-1}{s-1}\rfloor=2.$
Therefore~(\ref{que}) becomes:
\begin{equation}\label{final}
PoA(\rho, x_0^e, x_1^e)\ge 2 \cdot\Biggl(\frac{1}{\frac{1}{x^*}+\rho} \Biggr)\Biggl(1-\frac{1}{2}\Bigl(\frac{x^*}{x_e} \Bigr) \Biggr)\Biggl(\frac{1+\rho x_0^e+\rho^2 x_0^ex_1^e}{x_0^e+\rho (x_1^e)^2}\Biggr).
\end{equation}

Now, $x^*$ and $x_e,$ both relate to the loss system, and we have already shown that~(\ref{xstar0}) implies that $x_e\geq x^*,$ \,  and: $x^*\rightarrow \infty$  and  $\frac{x^*}{x_e}\rightarrow 0,$ \ when $x_e\rightarrow\infty.$
Note that when $s\to\infty,$ then $x_0^e\sim s^2\to\infty.$
Hence when $s$ goes to infinity, the first two factors in~(\ref{final}) converge to $\frac{1}{\rho}$ and $1,$ respectively, and so:
$$\lim_{s\to\infty}PoA(\rho, x_0^e, x_1^e)\ge \lim_{s\to\infty}
\frac{1+\rho x_0^e+\rho^2 x_0^ex_1^e}{\rho(x_0^e+\rho (x_1^e)^2)}=\lim_{s\to\infty} \frac{s^3}{s^2}= \infty,$$
which completes the proof.
\hfill\quad$\Box$

\section{Acknowledgement}
We thank the reviewers for their thorough review. We highly appreciate the comments and
suggestions, which significantly contributed to improving the quality of this publication.

\end{document}